\documentclass[runningheads]{llncs}

\usepackage{amsmath}
\usepackage{amsfonts}
\usepackage{amssymb}
\usepackage{hyperref}
\usepackage{graphicx}

\newcommand{\R}{\mathbb{R}}               
\newcommand{\M}{\mathcal M}              
\newcommand{\dd}{\mathrm{d}}
\newcommand{\ftf}{\omega}                
\newcommand{\ftfl}{\omega_{\mathcal L}}                
\newcommand{\fof}{\eta}                  
\newcommand{\defn}{\equiv}               
\newcommand{\hvf}[1]{\setbox0=\hbox{$#1$}%
  \ifdim\wd0>1em\widehat{#1}\else\hat{#1}\fi} 
\newcommand{\field}[1]{\Gamma(#1)}       
\newcommand{\liealgebra}[1]{\mathfrak{\MakeLowercase{#1}}} 
\newcommand{\G}{\mathcal G}              
\newcommand{\g}{\liealgebra G}           
\newcommand{\K}{\mathcal K}              
\renewcommand{\k}{\liealgebra K}         
\newcommand{\p}{\liealgebra P}           
\newcommand{\tang}[1]{\partial_{#1}}                
\newcommand{\GL}[1]{\mathop{\rm GL}(#1)} 
\newcommand{\metric}[2]{\langle #1, #2 \rangle} 
\newcommand{\Ad}{\mathop{\rm Ad}\nolimits} 
\newcommand{\ad}{\mathop{\rm ad}\nolimits} 
\newcommand{\End}[1]{\mathop{\rm End(#1)}} 
\newcommand{\tr}{\mathop{\rm tr}}          
\newcommand{\legtr}{\mathbb F \mathcal L}

\newcommand{\dt}{\delta t}                 
\newcommand{\refeq}[1]{(\ref{#1})}         
\newcommand{\refsec}[1]{\S\ref{#1}}        
\newcommand{\U}{V}                         

\def
\rational#1#2{{\mathchoice{\textstyle{#1\over#2}}%
  {\scriptstyle{#1\over#2}}{\scriptscriptstyle{#1\over#2}}{#1/#2}}}
\def\half{\rational12}		

\begin{document}
\title{Hamiltonian Monte Carlo On Lie Groups and Constrained Mechanics on Homogeneous Manifolds}
\titlerunning{HMC on Lie Groups and Homogeneous Manifolds}
%
\author{Alessandro Barp\inst{1,2}
}
\authorrunning{Barp Alessandro}
%
\institute{Imperial College London, Kensington, London SW7 2AZ, UK \and
The Alan Turing Institute, 96 Euston Rd, Kings Cross, London NW1 2DB UK 
\email{a.barp16@imperial.ac.uk}}
\maketitle              
\begin{abstract}
In this paper we show that the Hamiltonian Monte Carlo method for compact Lie groups constructed in \cite{kennedy88b} using a symplectic structure
can be recovered from canonical geometric mechanics with a bi-invariant metric.
Hence we obtain the correspondence between the various formulations of Hamiltonian mechanics on Lie groups, and their induced HMC algorithms.
Working on $\G\times \g$ we recover the Euler-Arnold formulation of geodesic motion, and construct explicit HMC schemes that extend \cite{kennedy88b,Kennedy:2012} to non-compact Lie groups by choosing metrics with appropriate invariances. 
Finally we explain how mechanics on homogeneous spaces can be formulated as a constrained system over their associated Lie groups, and how in some important cases the constraints can be naturally handled by the symmetries of the Hamiltonian.

\end{abstract}
\section{Introduction}

The HMC algorithm~\cite{Duane:1987} is a method that generate samples from a probability distribution which has a density known up to constant factor with respect to a reference measure~\cite{Betancourt:2017a,Barp:2018,Betancourt:2017b}.
The method  was originally extended to compact Lie groups in \cite{kennedy88b}, and Poisson backets, shadow Hamiltonians and higher order integrators were discussed in \cite{Kennedy:2012}. 
The algorithm was designed by  constructing a symplectic structure (and thus a mechanics) on $\G \times \g^*$,
but the relation with the standard mechanics induced by the Liouville form on $T^*\G$ was never explained.
One of the difficulties that arise when sampling distributions on manifolds using HMC is the requirement to compute the geodesic flow accurately in order to maintain the motion on the manifold \cite{byrne2013geodesic,holbrook2018geodesic,barp2019homogeneous}.

In the context of Lie groups, the geodesic motion was handled in  \cite{kennedy88b,Kennedy:2012,clark2008tuning,knechtli2017lattice} by assuming the Killing form defined a positive-definite Riemannian metric, in which case geodesics are given by 1-parameter subgroups.
This holds for $SO(n)$ and $SU(n)$ but is not case for most other Lie groups, such as $\GL n$ or $SL(n)$.\footnote{In~\cite{Kennedy:2012} the equations of motion are given for semisimple lie algebras, but the derivation uses the Killing form which is not Riemannian in general and as a result cannot be used for HMC (due to the momentum refreshment step).}

In this article we formalise the HMC schemes given in~\cite{Kennedy:2012,kennedy88b} and connect them to standard geometric mechanics on Lie groups, which allows us to remove the compactness assumption on the Lie algebra (or equivalently, the assumption that the Killing form is negative definite).

In \refsec{Symplectic-structures-groups} we show the symplectic structure used in \cite{Kennedy:2012,kennedy88b} is symplectomorphic to the canonical one on $T^*\G$. 
It follows that many results from geometric mechanics that assume standard mechanics, such as symplectic reduction, can be used in HMC~\cite{abraham1978foundations}.
Since HMC is usually implemented on the tangent bundle, we derive the corresponding symplectic structure on $\G \times \g$ induced by a metric. 
In \cite{Kennedy:2012,kennedy88b} the kinetic flow was simply given by left-invariant vector fields as a result of the bi-invariance of the Killing-form. 
This corresponds to solving the reconstruction equation for a constant curve on $\g$, 
but, as we show in \refsec{ham-flow-arnold}, the geodesic flow of more general left-invariant metrics
must solve both the Euler-Arnold and reconstruction equations.
Hence in \refsec{geodesics-groups} we consider the Euler-Arnold equation and explain how it can be simplified and solved by choosing inner products with appropriate symmetries.
In particular we derive explicit HMC schemes for $\GL n$ and semi-simple Lie groups.
In \refsec{leap-lie} we provide the explicit relation between the leapfrog scheme on Lie groups as used in \cite{barp2019homogeneous} and the general one on manifolds; 
while \refsec{force-integrator} discusses more efficient integrators. 
In \refsec{constrained-hom} we explain how more generally the mechanical system on homogeneous manifolds can be derived from a constrained mechanics defined on Lie groups.
For reductive homogeneous manifolds, the constraints can be naturally handled by choosing a Hamiltonian with sufficient symmetries. 
Hence we recover the HMC scheme proposed in~\cite{barp2019homogeneous} to sample distributions on naturally reductive homogeneous spaces using Lie group mechanics and symplectic reduction.  
We will see that this is still true for the discretised motion defined by the leapfrog method or force-gradient integrators (i.e., these preserve the symmetries).

\section{Hamiltonian Monte Carlo on Lie Groups}
 \label{HMC-lie-groups}

\subsection{Lie Groups}

Let $L_g :\G\to \G$ denote the left translation $L_gh=gh$ on a Lie group $\G$, and $\g$ (and $\g_L$) denote the lie algebra (of left-invariant vector fields) (see \refsec{lie-algebra} for more details). 
The Maurer-Cartan form $\theta \in \field{T^*\G \otimes \g}$ \footnote{$\field{\M}$ denotes the set of smooth sections of a bundle $\M$. If the base space $\mathcal B$ is ambiguous we write $\field{\mathcal B, \M}$. } 
is defined by $\theta:g \mapsto \tang g L_{g^{-1}}$ where $\tang{}$ denotes the tangent map. 
Let $\xi_i$ be a basis of $\g$, $e_i\in \g_L$ be the induced left-invariant vector fields (so $e_i(1)=\xi_i$), and $\theta^i$ the dual 1-forms, $\theta^i(e_j)=\delta^i_j$.
We will adopt the Einstein notation and sum over repeated indices. 
\subsection{Symplectic Structures on Lie Groups}\label{Symplectic-structures-groups}

In this section we derive the symplectic structures on $\G\times \g^*$, $T\G$ and $\G \times \g$ induced by the canonical 1-form on $T^*\G$.

We define the momentum as $p\defn \big(\theta^{-1}\big)^*:T^*\G \to \g^*$.
 In a basis $e_i$ of $\g_L$ we may expand these as
 $\theta = e_i(1) \otimes \theta^i$ and $p= \theta^i(1)\otimes e_i$ where we use the canonical isomorphism $T_g\G \cong T^{**}_g\G$ to view $e_i$ as a function $T^*\G \to \R$, i.e., $e_i: (g,\alpha_g) \mapsto \alpha_g \big(e_i(g)\big)$.
We can then define the 1-form $\fof \defn e_i \pi^* \theta^i \in \field{T^*\G, T^*T^*\G}$ on $T^*\G$ (here $\pi:T^*\G\to\G$), and it is easy to check it yields the canonical Liouville form~\cite{barp2019homogeneous}. 
The Liouville form may also be expressed in a frame-independent manner as $p\big( \pi^*\theta\big)$, as well as in terms of the right Maurer Cartan form~\cite{alekseevsky1994poisson}.
The symplectic structure is its exterior derivative\footnote{Here we use the fact that $\theta$ is a flat gauge-field, and thus satisfy Maurer-Cartan equation
$0=\dd \theta + \half [\theta \wedge \theta]= e_k(1) \otimes \big(\dd \theta^k+\half c^k_{ij}\theta^i\wedge \theta^j \big)  $. }
\begin{equation} \label{ftf-lie-group}
 \ftf \defn -\dd\fof = -\dd(e_i \pi^{*}\theta^i)
  = \pi^{*}\theta^i\wedge \dd e_i
    + \half  e_ic^i_{jk}\pi^{*}\theta^j\wedge\pi^{*}\theta^k,
\end{equation} 
and a similar expression holds for any symplectic fiber bundle $\pi:(E,\omega)\to \M$ over a parallelisable manifold $\M$ (i.e., manifolds for which $T\M$ is trivial)~\cite{alekseevsky1994poisson}.
 The symplectic volume form is then given by the product of the Haar measure $\theta^1 \wedge \cdots \wedge \theta^n$ on $\G$ with the Lebesgue measure on the fibres
  $$\omega^n \propto \dd e_1 \wedge \cdots \wedge \dd e_n  \wedge\pi^*\theta^1\wedge \cdots \wedge \pi^*\theta^n $$

Originally HMC on compact matrix groups was constructed on the trivial bundle $\G\times \g^*$, and although the momentum coordinates were left undefined,\footnote{Technically in \cite{Kennedy:2012} the symplectic structure is defined as $\dd p$ for $p\in T_g^*\G$ which does not define a differential form over $T^*\G$  since the components $p_i$ are then constant.} we can recover the original mechanical system by defining
 $\tilde p_i : \G \times \g^* \to \R$ s.t. $\tilde p_i : (g,\alpha_1) \to \alpha_1(\xi_i)$, from which we can form the symplectic structure $\tilde \omega = -\dd(\tilde p_i \pi^* \theta^i)$.
Clearly $\omega$ and $\tilde \omega$ are symplectomorphic with respect to the canonical isomorphism $T^*\G \cong  \G \times \g^*$, $(g,\alpha_g)\mapsto (g, L^*_g \alpha_g)$ (see \refsec{symplecto-cot}).

In practice however we usually work on the tangent bundle. 
Consider a Riemannian metric $\metric{\cdot}{\cdot}$ on $T\G$. 
The associated musical isomorphism, $\legtr:T\G \to T^*\G$, $v \mapsto \metric{v}{\cdot}$, enable us to pull-back $\omega$ to a symplectic structure on $T\G$.
 
The symplectic structure $\ftfl$ induced on $T\G$ by the musical isomorphism is 
the exterior derivative of the 1-form $-v^j\pi^*\big(\metric{e_i}{e_j}\theta^i\big)\in \field{T\G,T^*T\G}$, where $v^j:T\G \to \R$ satisfies $(g,u_g) \mapsto \theta^i_g(u_g)$, and now $\pi:T\G \to \G$ is the tangent bundle projection (see \refsec{pull-tan}). 
To take advantage of the symmetries of $\G$, we assume $\metric{\cdot}{\cdot}$ is left-invariant and set $g_{ij}\defn\metric{\xi_i}{\xi_j}_1  $.
Then
\begin{equation} \label{ftf-lie-group}
 \ftfl \defn \legtr^*\ftf = -\dd(g_{ij}v^j \pi^{*}\theta^i)
  = g_{ij}\pi^{*}\theta^i\wedge \dd v^j
    + \half  g_{ir}v^rc^i_{jk}\pi^{*}\theta^j\wedge\pi^{*}\theta^k.
\end{equation}

As above, $\ftfl$ can be pulled-back to define an equivalent mechanics on $\G\times \g$ 
(the coordinate functions $v^j$ are now the maps $(g,X)\mapsto \theta^j_1(X)$).
In particular, if the metric is left-invariant the Hamiltonian $H\defn V \circ \pi + \half  \metric{\cdot}{\cdot} $ on $T\G$ defines an equivalent Hamiltonian $V \circ \pi + \half  \metric{\cdot}{\cdot}_1$ on $\G \times \g$ (see \refsec{symplecto-tan}).

\subsection{Hamiltonian Fields and Euler-Arnold Equation} \label{ham-flow-arnold}

We now focus on the mechanical system on $\G\times \g$ induced by a left-invariant metric $T=v^iv^jg_{ij}:\g\to \R$.
For any function $f:\G\times \g \to \R$, its Hamiltonian vector field $\hvf f \in \field{T\G \oplus T\g}$ is defined by $\ftfl(\hvf f,\cdot)=\dd f$.
We have (see  \refsec{proof-hamvf})
\begin{equation}
\hvf f = g^{jk} \xi_k(f) e_j + \Big( g^{ik}g^{lj} v^r c_{rji} \xi_l(f)-g^{jk}e_j(f) \Big) \xi_k.
\end{equation} \label{hvf-f}
For the potential energy and kinetic energy induced by a left-invariant metric
 $$ \hvf T =  v^k e_k+g^{ik} v^rv^s c_{rsi} \xi_k \qquad \hvf V = - g^{jk}e_j(V) \xi_k$$

The Poisson brackets associated to $\ftfl$ are defined by
 $\{r,h\} \defn  \ftfl (\hvf r, \hvf h)$ for any functions $r,h:\G\times \g\to \R$, 
 and describe
the rate of change of $r$ along the flow of $h$ since
$\{r,h\} =  \ftfl (\hvf r, \hvf h)=\dd r(\hvf h)= \hvf h[r]$. 
This is often written $\dot r = \{r,h\}$ where $\dot r:\G\times \g \to \R$ is the function giving the rate of change of $r$  in the direction $\hvf h$,
 $\dot r(u)\defn \frac{\dd r(\Phi^{\hvf h}_u(t))}{\dd t}|_{t=0}$ (here $\Phi^{\hvf h}:\R\times (\G \times \g) \to \G\times \g$ denotes the flow of $\hvf h$).
It follows immediately that a function does not vary along its Hamiltonian flow: $ \dot h = \{h,h\}=0$, and in fact neither does the symplectic structure $\Phi^{\hvf h}(t)^*\ftfl =\ftfl$, 
which enables the construction of symplectic integrators (see \refsec{Symplectic-integrator}). 
It is also common to write the Poisson bracket of a curve $\gamma:[a,b] \to T\G$ with a function $r$ by setting 
$\dot \gamma = \{\gamma, h\} \defn \hvf h( \gamma)$.
This is a differential equation stating that $\dot \gamma$ is the integral curve of $\hvf h$. 
In particular, a curve satisfying $\dot \gamma = \{ \gamma,H\}$ with initial conditions $(g_0,v_0)=\gamma(0)$ is known as a solution of Hamilton's equations.

An integral curve $\gamma_V$  of $\hvf V$ satisfies $\dot \gamma_V = \hvf V (\gamma)$. 
Let $\gamma_V(t) =\big(q(t),v(t)\big)$.
The equation then reads $\big(\dot q, \dot v \big) = \big(0,-g^{jk}e_j(V) \xi_k \big)\big|_{\big(q(t),v(t)\big)} $, which implies 
\begin{equation}
q(t) =q(0) \qquad  v(t) =v(0)-tg^{jk}e_j(V)|_{q(0)} \xi_k .
\end{equation}\label{vel-update-lie}
For the integral curves of $\hvf T$, 
we have $\big(\dot q, \dot v \big) =\big(v^k e_k, g^{ik} v^rv^s c_{rsi} \xi_k\big)|_{\gamma_T(t)}$.
The mysterious term $g^{ik} v^rv^s c_{rsi} \xi_k$ disappeared in both \cite{Kennedy:2012,kennedy88b} (and didn't appear in \cite{knechtli2017lattice}) since the structure constants of the Killing form are totally anti-symmetric (see \refsec{AdG}).
It is in fact nothing else than the $\metric{\cdot}{\cdot}$-adjoint of the adjoint representation.
Recall that the adjoint representation $\ad:\g \to \End \g$ satisfies $\ad_{\xi} \defn \ad(\xi) = [\xi,\cdot]$. 
For each $\xi \in \g$ we can define its adjoint $\ad^*_{\xi}$ with respect to the inner product 
$\metric{\ad^*_{\xi}X}{Y} = \metric{X}{\ad_{\xi}Y}$. 
It follows that $\ad^*$ is bilinear so 
\begin{align*}
\metric{\ad^*_vv}{\xi_k} &= v^i v^j \metric{\ad^*_{\xi_i}\xi_j}{\xi_k}= v^i v^j \metric{\xi_j}{\ad_{\xi_i}\xi_k}
=
v^i v^j \metric{\xi_j}{[\xi_i,\xi_k]}\\
&=
v^i v^j \metric{\xi_j}{c_{ik}^r\xi_r}= c^r_{ik}v^iv^j g_{jr}= v^i v^jc_{jik}.
\end{align*}

Moreover $\metric{\ad^*_vv}{\xi_k}=\metric{\big(\ad^*_vv\big)^s\xi_s}{\xi_k}
=\big(\ad^*_vv\big)^s g_{sk}$. 
Combining these two equations gives 
$\ad^*_vv= g^{kl} v^i v^jc_{jik} \xi_l$.
Hence the integral curves of $\hvf T$ satisfy 
 $$\big(\dot q, \dot v \big) =\big(v^k e_k, \ad^*_vv\big)|_{\gamma_T(t)}.$$ 
 
 The equation $\dot v = \ad^*_vv$ is a first order differential equation on the lie algebra $\g$ known as the Euler-Arnold equation.
 In order to derive the integral curves of $\hvf T$ we can proceed in two steps.
 We first solve the Euler-Arnold equation to give a solution $s:[0,T] \to \g$. 
 The equation for $\dot q$ then states that $\dot q =v^ke_k|_{(q(t),s(t))}= s^k(t)e_k(q(t)) = \tang 1 L_{q(t)} s(t)$, which is a first order ODE for $q$ called the reconstruction equation.

In summary we have shown that the Hamilton's equation $\dot \gamma = \hvf H(\gamma)$ defined by the symplectic structure on $\G\times \g$ reads 
$$  \dot q= v^je_j(q)= \tang 1 L_q v \qquad \dot v = \ad^*_vv- g^{jk}e_j(V)(q) \xi_k,  $$
which can be recognised as the equations arising from the Hamilton-Pontryagin variational principle.

We finally note that the left-invariance of the metric enabled us to obtain two first order differential equations for the geodesic flow, rather than the usual second-order differential equation on a general Riemannian manifold.
This is known as the Euler-Poincaré reduction.

\subsection{Leapfrog Integrator} \label{leap-lie}

On an arbitrary Riemannian manifold $\M$, we can define the leapfrog integrator for the Hamiltonian system $H=V+T$ as a composition of the kinetic $\hat T$ and potential $\hat V$ fields defined by the canonical symplectic structure on $T\M$~\cite{Leimkuhler:1996}.
The potential flow leads to the velocity update
\begin{equation}
q(t) =q(0) \qquad v(t)=v(0)-t\nabla_{q(0)} V ,
\end{equation}\label{velo-update}
(here $\nabla$ is the Riemannian gradient) and the kinetic (or geodesic) flow leads to the update
$$ \big(q(t),v(t)\big) = \big( \exp_{q(0)}(v(0)), \mathbb P_{q(t)}v(0) \big),$$
where $\exp_qv$ denotes the Riemann exponential map and $\mathbb P_q v$ the parallel transport along the geodesic.
Let us see the explicit relation between these updates and the ones derived in \refsec{ham-flow-arnold}.
Suppose that $\M=\G$ is a Lie group and $T$ is defined by a left-invariant metric $\metric{\cdot}{\cdot}$. 
From \refsec{Symplectic-structures-groups} we know that, under the symplectomorphism $T\G \cong \G \times \g$, the Hamiltonian vector-field induced by the inner product $\metric{\cdot}{\cdot}_1:\g\to \R$ push-forwards to the geodesic flow induced by $T:T\G\to \R$. 
More precisely if $(q(t),v(t)) \in \G \times \g$ is the integral curve of the Hamiltonian vector field of $\metric{\cdot}{\cdot}_1$, then the geodesic flow of the Riemannian metric $T=\metric{\cdot}{\cdot}$ on $T\G$ is given by $\big(q(t),\tang 1 L_{q(t)}v(t)\big)= \big(q(t), \dot q(t) \big)$.

For the potential term, observe that\footnote{Note that if $\alpha \in \field{T^*\G}$, then $\alpha = e_i(\alpha) \theta^i= \alpha(e_i) \theta^i = \alpha_i \theta^i$. 
Moreover $u = g^{ij}\alpha_ie_j$ is the vector field associated to $\alpha$ by the metric since 
$\alpha(v) = \alpha_i \theta^i(v^je_j)= \alpha_i v^i = \alpha_k \delta^k_i v^i = \alpha_k g^{kr}g_{ri} v^i = u^rg_{ri}v^i = \metric uv $. }
 $$\nabla_gV= g^{jk}(g) \dd_g V(e_k) e_j(g)=g^{jk}(1)e_k(V)(g) \tang 1 L_g \xi_j$$
  (by left-invariance $g^{ij}$ is constant), and $v(t)=\theta^i_g(v(t))e_i(g)=v^i(t) \tang 1 L_g \xi_i$.
Thus under the symplectomorphism $ \G \times \g \to T\G $, we have ($\xi_i \to \tang 1 L_g \xi_i )$ and \refeq{velo-update} push-forwards to the velocity update  \refeq{vel-update-lie}.

Finally let us see how the velocity update
\refeq{vel-update-lie} is computed in practice when $\G$ is a matrix group.
In that case we need an extension $W$ of $\U$ defined on an open subset of $\GL n$. 
Then
$$e_i(V)(g) = 
 \dd(W\circ \iota)(  g\cdot \xi_i) =\frac{\partial W}{\partial {x^{jk}}}\big(\tang{}( \iota\circ L_g ) \xi_i)\big)_{jk} = \tr\Big( \frac{\partial W}{\partial {x}}^T\cdot g\cdot \xi_i\Big),$$
 where $\cdot$ is matrix multiplication, $\iota:\G \to \GL n$ is the inclusion (or a more general injective representation), and in the last line we have identified $ g \sim \iota(g), \xi_i \sim \tang{}\iota ( \xi_i)$.
 We arrive at the general velocity update formula for a matrix group
$$v^i(t)\xi_i=v^i(0) \xi_i-tg^{ij}(1)\tr\Big( \frac{\partial W}{\partial {x}}^T\cdot g\cdot \xi_j\Big)\xi_i.$$ 

In lattice gauge theory we typically have $W(x)= Re \tr\big(Ux\big)$ for a constant complex matrix $U$.

\subsection{Efficient Integrators} \label{force-integrator}

While the leapfrog integrator is the most common choice of integrator, it is not always the optimal choice.
For example \cite{takaishi2006testing,omelyan2001algorithm} consider a symplectic integrator of the form\footnote{Using the notation $\exp(t \hvf H) \defn \Phi^{\hvf H}(t)$ for the Hamiltonian flow}
 $\exp \big(\lambda \dt \hvf T\big)\exp \big(\half \dt \hvf V\big)\exp \big( (1-2\lambda) \dt \hvf T\big)\exp \big(\half \dt \hvf V\big)\exp \big(\lambda\dt \hvf T\big)$ which, despite being computationally more expensive, they show is $50\%$ more efficient than leapfrog.
In \cite{Kennedy:2012} a ``force-gradient" integrator is examined which involves the Hamiltonian flow $\hvf{\{V, \{V,T\} \}}$, and thus contain second-order information (derivatives) about $V$. 
The question of tuning HMC using Poisson brackets to minimise the cost was considered in \cite{clark2008tuning}; see also \cite{celledoni2014introduction} for a discussion on Lie group integrators, and~\cite{bou2018geometric} for integrators which outperform leapfrog on $\R^n$.

\section{Geodesics} \label{geodesics-groups}
We now analyse how Euler-Arnold equation $\dot v = \ad^*_vv$ simplifies when the inner product on $\g$ is $\ad_{\k}$-invariant for a subalgebra $\k\subset \g$.
When $\G$ is connected, $\ad_\k$-invariant inner products correspond to left $\G$-invariant and right $\K$-invariant Riemannian metrics on $\G$. \footnote{i.e., the left and right actions $L_g$ and $R_k$ are isometries for $g\in\G,k\in \K$ (here $\K$ is the Lie group associated to $\k$).
When $\G$ is not connected an analogous correspondence exists with $\Ad_{\K}$-invariant inner products~\cite{barp2019homogeneous}, i.e., inner products s.t. $\metric{\Ad_k u}{\Ad_k v}_1 =\metric uv_1$ for $u,v \in \g$, which are automatically $\ad_\k$-invariant.}
When $\k = \g$ these correspond to bi-invariant Riemannian metric on $\G$, which always exists when $\G$ is compact~\cite{Nielsen:2012:MIG:2371218}.

\subsubsection{$\Ad \G$-invariant Inner Products} \label{AdG} The simplest case arises when $\ad^*_vv=0$. 
This situation occurs for example when the inner product on $\g$ is $\ad$-invariant, 
$\metric{\ad_u v}{w}+\metric{v}{\ad_uw}=0$, which implies $\ad^*=-\ad$ and $\ad^*_vv=-[v,v]=0$. 
We can find $\ad$-invariant inner products on $\g$ when $\G$ has an $\Ad \G$-invariant inner product, or equivalently if it has a bi-invariant metric.
A particular $\Ad \G$-invariant inner product is available on Lie groups whose lie algebra is compact, which means the Killing-form $B(u,v) \defn \tr\big( \ad(u)\ad(v) \big)$ is negative definite, and so $-B$ defines a positive definite $\ad $-invariant inner product. 

When $\ad^*_vv=0$, then $v(t)=v_0$, so the reconstruction equation becomes $\dot q =  v^ke_k|_{(q(t),v_0)} = v_0^k e_k(q(t))= \tang 1 L_{q(t)} v_0 = \frac{\dd}{\dd t} L_{q(t)}e^{tv_0} |_{t=0} $, which implies $q(t) = q(0)e^{tv_0}$ (here we used the fact that $e^{tv_0}$ is tangent to $v_0$ at $t=0$).
Hence we obtain the kinetic flow as in \cite{Kennedy:2012} $$\gamma_T(t) = \big(q(0)e^{tv_0},v_0).$$

\subsubsection{Reductive Decomposition}

Sometimes it is not possible to find a bi-invariant metric on $\G$. 
The natural next step is to consider inner products on $\g$ which are $\ad_\k$-invariant, where $\k$ is a subalgebra of $\g$.
This means that $\metric{\ad_u v}{w}+\metric{v}{\ad_uw}=0$ for all $u \in \k$, and it follows that $\ad^*_u = -\ad_u$.
Then $\K$ is totally geodesic in $\G$ and its geodesics are given by 1-parameter subgroups (see \refsec{tot-geo}). 
Let $\p \defn \k^{\perp}$, so $\g =\k \oplus \p$ and $v=v_\k+v_\p$.
The Euler-Arnold equation reads
$\dot v = \dot v_\k + \dot v_\p = \ad^*_vv=\ad^*_{v_\k}v_\k+\ad^*_{v_\k}v_\p+\ad^*_{v_\p}v_\k+\ad^*_{v_\p}v_\p=-[v_\k,v_\p]+\ad^*_{v_\p}v_\k+\ad^*_{v_\p}v_\p$. 
We can simplify this further if the decomposition $\k \oplus \p$ is naturally reductive, i.e., $[\k,\p]\subset \p$ and for all $u,v,w\in\p$,  $\metric{[u,v]_\p}{w}= \metric{u}{[v,w]_\p}$~\cite{alekseevsky2007riemannian}. 
This implies $\ad^*_{v_\p}v_\p =0$ (see \refsec{natural-geo}), and we are left with $$\dot v =-[v_\k,v_\p]+\ad^*_{v_\p}v_\k.$$
   Notably the constant curve $v=v_0$ with either $v_0 \in \k$ or $v_0 \in \p$ is a solution.

\subsubsection{Matrix Groups} \label{mat-groups}

Consider a subgroup $\G \subset \GL n$ and set $\g = \p \oplus \k$ where $\p$ is a subspace of the space of symmetric matrices, while $\k = \mathfrak{so}(n)$.
We have in mind the cases $\GL n, SL(n)$ and $SO(n)$ with $\p$ the space of symmetric matrices, the space of traceless symmetric matrices, and $\{0\}$ respectively.
We equip $\g$ with the standard inner product $\metric{A}{B} \defn \tr(A^TB)$ which is $\Ad \K$-invariant. 
Since the product of a symmetric and an antisymmetric matrix is traceless, we have
$\metric{A}{B}= \tr(A_\p B_\p)-\tr(A_\k B_\k)$.

Note that for $\G=SO(n)$, this is just the negative of the Killing form (up to a positive constant);
while on $\GL n$ and $SL(n)$, the Killing form, respectively given by $2n\tr(AB)-2\tr(A)\tr(B)$ and $2n\tr(AB)$, doesn't define an inner product (it is degenerate and pseudo-Riemannian respectively~\cite{miolane2015computing}).

For our inner product we have $\ad^*_{S}A=[S,A]$ for any $S\in\p, A\in \k$, hence Euler-Arnold equation simplifies to $\dot v = -2[v_\k,v_\p]$ with solution (see \refsec{geo-matrix}) 
\begin{equation}
\big(q(t),v(t)\big)=\Big(q(0)e^{\big(v_\p(0)-v_\k(0)\big)t}e^{2v_\k(0) t}, v_\k(0)+ \Ad e^{-2 v_\k(0)t}(v_\p(0)) \Big).
\end{equation}

\subsubsection{Semi-simple Symmetric Spaces}

The equation for the geodesics above also holds  for a semi-simple Lie group $\G$ using an inner product based on the Killing form $B$ (note $\GL n$ is not semi-simple). 
Let us consider a subgroup $\K$ for which $\G/\K$ is a symmetric space. 
In that case we can find a Cartan decomposition $\g=\k\oplus \p$.
By corollary 5.4.3 \cite{jost2008riemannian} the Killing form is negative-definite on $\K$.
We say $\G/\K$ is of compact (non-compact) type if $B$ is negative (positive) definite on $\p$.
If $\G/\K$ is of compact type the inner product $\metric{v}{u} \defn B(v_\k,u_\k)+B(v_\p, u_\p)$ is $\Ad \G$-invariant, while if it is of non-compact type  $\metric{v}{u} \defn -B(v_\k,u_\k)+B(v_\p, u_\p)$ is $\Ad \K$-invariant.
Then $\ad^*_{v_\p}v_\k=[v_\p,v_\k]$ (see Lemma 5.5.4 \cite{jost2008riemannian}) and we can find the integral curves of $\hvf T$ as in \refsec{mat-groups}.

\subsubsection{Approximating the Matrix Exponential}

In practice there are several ways to approximate the matrix exponential, for example by combining a Padé approximant with a projection from $\GL n$ to $\G$~\cite{bou2009hamilton}.
In particular the Cayley transform $z \mapsto (1+\half z)/(1-\half z)=(1-\half z)^{-1}(1+\half z)$, which is the Padé-$(1,1)$ approximant, maps exactly onto $\G$ for quadratic groups\footnote{By a quadratic group (sometimes called $J$-orthogonal) we mean a subgroup $\G=\{ M \in \GL n : M^TJM=J\}$ for some $J\in\GL n$}(in fact diagonal Padé approximants are product of Cayley transforms)~\cite{iserles2001cayley} .

\section{Constrained Mechanics on Homogeneous Manifolds} \label{constrained-hom}
It was shown in \cite{barp2019homogeneous} that HMC on naturally reductive manifolds can be implemented using mechanics on Lie groups.
Consider a Homogeneous manifold $\G/\K$ (where $\K$ is a closed subgroup of $\G$) and a Hamiltonian $\tilde H=\tilde V+\tilde T$ on $T \big(\G/\K\big)$. Let $\pi:\G \to \G/\K$ be the quotient projection.
We can lift $\tilde H$ to $T\G$ by defining $H \defn (\tang{}\pi)^*\tilde H$.
While $V \defn \pi^* \tilde V$ is well-defined on $\G$, the lift $T \defn (\tang{}\pi)^* \tilde T$ defines a degenerate kinetic energy on $T\G$, since it maps all vectors in the vertical space $\text{ver}_g \defn \text{ker}(\tang{g}\pi)$ to zero.
The vertical space is described by the action fields $\text{ver}_g =\{\xi_{\G}(g) :\xi \in \k\}$ where $\xi_{\G}(g) = \frac{d}{dt}e^{t\xi}\cdot g |_{t=0} = \tang 1 L_g(\xi)$.

In order to circumvent the issues related to the degeneracy of $T$, we need to make a choice of horizontal space complementary to $\text{ver}_g$, or equivalently choose a connection on the principal bundle $\G \to \G/\K$. 
This defines a unique way to lift a curve on $\G/\K$ to a curve on $\G$, and a complementary subspace $\p$ to $\k$, i.e., $\g=\k\oplus \p$ (see \refsec{horizontal}). 
Then, if $(q,v)$ satisfies the constrained equation
$$  \dot q= v^je_j(q) \qquad \dot v = \ad^*_vv- g^{jk}e_j(V)(q) \xi_k \qquad v(t) \in \p, $$
the projected curve $\pi(q(t))$ solves Hamilton's equations on $\G/\K$~\cite{lee2017global}. 
We note the constraint $v(t) \in \p$ can be rewritten in the form $ver\big( v(t)\big)=0$ where $ver:\g \to \k$ is the vertical projection.

If $\G/\K$ is a reductive homogeneous manifold, there is a one-to-one correspondance between $\Ad_{\G}\K$-invariant inner products on $\p$ and $\G$-invariant metric on $\G/\K$. 
In particular if we extend the inner product on $\p$ to a non-degenerate quadratic form on $\g$ s.t. $\p = \k^{\perp}$, it defines a pseudo-Riemannian metric on $\G$ that is left $\G$-invariant and right $\K$-invariant,
 and a kinetic energy $T$ with $T = (\tang{}\pi)^* \tilde T$ on $\p$.
The resulting Hamiltonian $ H=V+T$ is right $\K$-invariant so by Noether's theorem, if $v(0) \in \p$, then $v(t) \in \p$ for all $t$. 
This means the constraint $v(t)\in \p$ is naturally handled by the symmetries of $H$ (see \refsec{sym-Ham}).

If $\G/\K$ is further naturally reductive, since $\ad^*_vv=0$ for $v\in \p$ the above system becomes 
$$  \dot q= v^je_j(q) \qquad \dot v =- g^{jk}e_j(V)(q) \xi_k \qquad v(0) \in \p, $$
and we recover the HMC algorithm on naturally reductive homogeneous spaces proposed in \cite{barp2019homogeneous}.

It is clear that integrators that are built by alternating between the flows of $V$ and $T$ are automatically horizontal (the trajectory stays in $\G \times \p$)  since both $V$ and $T$ are $\K$-invariant.
This is still true for the force-gradient integrators \refsec{force-integrator}.
Indeed, from the $\ad_\k$-invariance of the inner product on $\g$, we have $c_{rij}=-c_{jir}$ for any $i \in \k$ (the structure constants always satisfy $c^i_{jk}=-c^i_{kj}$).
Hence if $V$ is invariant under the right action of $\K$ on $\G$, then $\hvf{\{V, \{V,T\} \}}= -2g^{jk}g^{ls} e_l(V) e_je_s(V)  \xi_k\in \p$ (see \refsec{Force-proof}).

\section{Acknowledgements}  
We thank Prof. Anthony Kennedy and Prof. Mark Girolami for useful discussions.
AB was supported by a Roth Scholarship funded by the Department of Mathematics, Imperial College London, by EPSRC Fellowship (EP/J016934/3) and by The Alan Turing Institute under the EPSRC grant [EP/N510129/1].

%
%
%
\bibliographystyle{splncs04}
\bibliography{GSIref}

\section{Supplementary Material}
\subsection{Symplectomorphism $T^*\G \cong \G \times \mathfrak{g}^*$} \label{symplecto-cot}

The cotangent bundle over $\G$ is  isomorphic to a trivial bundle 
$T^*\G \cong \G \times \g^* \cong \G \times \R^n$. Indeed fix a frame $e_i \in \g_L$ and co-frame $\theta^i \in \g^*_L$. We have the maps 
$$ (g, p) = (g,p_i \theta^i_g) \mapsto (g, p_i \theta^i_1) \mapsto \big(g,(p_1,\ldots,p_n) \big)$$

Letting $F:T^*\G \to \G \times\g^*$ be this diffeomorphism. Then 
$$ \big(F^* \tilde p_i \big) (g,\alpha_g) =(\tilde p_i \circ F)  (g,\alpha_g) =
\tilde p_i \circ ( F (g,\alpha_g) )=
\tilde p_i \circ (g, L^*_g\alpha_g)
= 
 L^*_g\alpha_g(e_i(1)) 
  = p_i(\alpha_g)$$
and $$F^*(\pi^* \theta^i)= \theta^i \circ \tang{}\pi \circ \tang{}F = \theta^i \circ \tang{}(\pi \circ F)= \theta^i \circ \tang{} \pi =  \pi^*\theta^i.$$

Hence $T^*\G$ with canonical symplectic structure is symplectomorphic to $\G\times \g^*$ with the symplectic structure from \cite{Kennedy:2012,kennedy88b}.

\subsection{Symplectic Structure on $T\G$} \label{pull-tan}

The pull-back of the coordinate functions $e_i$ is $e_i \circ \legtr: T\G \to \R$, and 
\begin{align*}
e_i \circ \legtr(h,v_h)&= e_i \big( \metric{v_h}{\cdot}\big)=\metric{v_h}{e_i(h)} = v_h^j \metric{e_j(h)}{e_
i(h)}\\
&= v^j_h \metric{e_j}{e_
i}(h)=\metric{e_j}{e_
i}(h) \theta^j_h (v_h)=\big(\theta^j \pi_{T\G}^*\metric{e_j}{e_
i} \big) (h,v_h),
\end{align*}
where $\pi_{T\G}:T\G\to \G$ is the canonical projection. 
Note in the final expression $\theta^j$ is viewed as a map $T\G \to \R$ (which we denote by $v^j$) rather than an element of $\field{T^*\G}$.
This gives the definition of the coordinate functions $v^j$ on $T\G$.
On the other hand, $\legtr^* \pi^*\theta^i= (\pi \circ \legtr)^* \theta^i$ and $\pi \circ \legtr:T\G \to \G$ is simply $\pi_{T\G}$.

In summary, $$\legtr^*\big( e_i \pi^*\theta^i\big)= \big(e_i \circ \legtr\big) \legtr^*\pi^*\theta^i=
\big(v^j \pi_{T\G}^*\metric{e_j}{e_
i} \big) \pi^*_{T\G}\theta^i=v^j \pi^*_{T\G}\big(\metric{e_j}{e_
i} \theta^i\big)$$

\subsection{Symplectomorphism $T\G \cong \G \times \g$} \label{symplecto-tan}
The isomorphism $T\G \cong \G \times \g$  can be written as $(g,v) \mapsto (g,\theta_g(v)) =(g, \tang g L_{g^{-1}}v)$, where $\theta$ is the (left) Maurer-Cartan form. 
We can thus pull-back  our Hamiltonian system from $T\G$ to $\G \times \g$: let $f: \G \times \g \to T\G$ be this map, and let $E = V\circ \pi +\half  \metric{\cdot}{\cdot}$ be an energy function on $T\G$.
Then
$$f^*E= E \circ f = V\circ \pi \circ f+\half  \metric{\cdot}{\cdot} \circ f= V\circ \pi + \half\metric{f(\cdot)}{f(\cdot)}$$
and 
\begin{align*} \metric{f(g,v^i e_1(1))}{f(g,v^j e_j(1))} & = \metric{(g,v^i e_i(g))}{(g,v^j e_j(g))}=  \metric{v^ie_i(g)}{u^je_j(g)}_g \\
&= L^*_g\metric{v^ie_i(1)}{u^je_j(1)}  
\end{align*}
Thus if the metric is left-invariant
$$ f^*E = V\circ \pi + \half \metric{\cdot}{\cdot}_1.$$

\subsection{Hamiltonian Vector Fields} \label{proof-hamvf}

We proceed similarly to \cite{Kennedy:2012} but work on $\G\times \g$ which will enable us to identify the coordinate vector fields $\partial_{v^i}$ with $\xi_i$. An arbitrary vector field $Y$ on $\G\times \g$ may be expanded as $Y=Y^i e_i + \bar Y^i \partial_{v^i}$, where $v^i:(g, u_1) \mapsto \theta^i_1(u_1)$ are the fibre coordinates.
This follows because any vector tangent to $\G \times \g$ is given by the derivative of a curve $\gamma: [a,b] \to \G \times \g$ which may be written as $\gamma(t) = (q(t),u(t))$; 
then $\dot \gamma = (\dot q, \dot u)$, where $\dot q \in T\G$ and $\dot u \in T\g$ can be expanded in terms of the (global trivialising) frames $e_i$ and $\partial_{v^i}$.

The equation $\ftfl(\hvf f, Y) = \dd f(Y)$ must hold for every $Y$ if $\hvf f$ is the Hamiltonian vector field of $f$.
Note $$\pi^* \theta^i \wedge \dd v^j=
\pi^* \theta^i \otimes \dd v^j - \dd v^j \otimes \pi^* \theta^i$$ and set
 $\hat f = X^i e_i + \bar X^i \partial_{v^i}$ 
Then $\pi^* \theta^i \wedge \dd v^j(\hat f, Y)= X^i \bar Y^j - Y^i \bar X^j $, 
moreover $\pi^* \theta^i \wedge \pi^* \theta^k(\hat f, Y)= X^jY^k-Y^jX^k$, and finally $\dd f(Y) =  Y(f)=Y^ie_i(f)+ \bar Y^i \partial_{v^i}f$.
Hence $\ftfl(\hvf f, Y) = \dd f(Y)$ is equivalent to 
$$ g_{ij}\big(X^i \bar Y^j - Y^i \bar X^j\big) + \half  v^r c_{rjk}\big(  X^jY^k-Y^jX^k\big)=Y^ie_i(f)+ \bar Y^i \partial_{v^i}f$$
where $c_{rjk} \defn g_{rl}c^l_{jk}$. 
As this must hold for every $Y$, equating the coefficients of $Y^i$ and $\bar Y^i$ yields $g_{ij}X^i= \partial_{v^j}f$ and 
$g_{ij}\bar X^j  = v^rc_{rji}X^j- e_i(f)$ respectively, where we have used $c_{rjk} = -c_{rkj}$ since the Lie bracket of vector field is antisymmetric: $[Z,W]=-[W,Z]$.

It follows that 
$$\hat f = g^{jk} \partial_{v^k}f e_j + \big( g^{ik}g^{lj} v^r c_{rji} \partial_{v^l}f-g^{jk}e_j(f) \big) \partial_{v^k}. $$
 In particular if $T= \half  g_{ij}v^i v^j$ and $V$ is independent of fibre coordinates $v^i$, then 
 $$ \hvf T =  v^k e_k+g^{ik} v^rv^s c_{rsi} \partial_{v^k} \qquad \hvf V = - g^{jk}e_j(V) \partial_{v^k}$$
 where we have used $e_i[T]=0$ since $e_i[T]= \half  e_i\big( g_{jk}(g)\big)v^jv^k=
 \half  e_i\big( g_{jk}(1)\big)v^jv^k=0$ by left-invariance of the metric. 
 
 In fact the coordinate vectors $\partial_{v^k}$ are simply the generators,
 $\partial_{v^k} = \xi_k$, since both of them have the integral curve $\gamma:t\mapsto t\xi_k$.
 Indeed, if $v=(v^1,\ldots,v^n):\g \to \R^n$ denotes the coordinate chart
 $\dot \gamma[f] = \frac{\partial (f \circ v^{-1})}{\partial v^i} \frac{\partial v^i\circ \gamma}{dt}|_{t=0}= \frac{\partial (f \circ v^{-1})}{\partial v^k} = \frac{\partial}{\partial v^k} (f)$. 
 
 Note that for a general $f$, $e_j(f)$ is a function on $T\G$, however since $V$ does not depend on the fibre coordinates, $e_j(V)$ is just a function of $\G$.

\subsection{Totally Geodesic Subgroups}
\label{tot-geo}
Recall we call a subgroup $\K \subset \G$ totally geodesic, if the geodesics of $\K$ are geodesics of $\G$.

Suppose we have an $\ad_\k$-invariant inner product on $\g$. 
Let $\p \defn \k^{\perp}$, then for any $u \in \k, v \in \p$, $\metric{u}{\ad_u v}= -\metric{\ad_uu}{v}=-\metric{[u,u]}{v}=0$, so Theorem 2 \cite{modin2010geodesics} is satisfied.
In fact $\k$ is ``easy totally geodesic" in $\g$, i.e.,  $\ad_\k \p \subset \p$, since for any $u,w \in \k$, $v\in \p$ we find
$\metric{\ad_uv}{w}= \metric{v}{\ad^*_uw}=-\metric{v}{[w,u]}=0$ since $[w,u]\in \k$ and $v\in \p$, hence $\ad_u v \in \p$ (see def 4.1 \cite{modin2010geodesics}). 

Finally we note that the restriction of the inner product to $\k$ is $\ad$-invariant and thus the geodesics are given by the 1-parameter subgroups~\cite{modin2010geodesics}.

\subsection{Naturally Reductive Geodesics}
\label{natural-geo}
Recall that we have an $\ad_\k$-invariant inner product, i.e. $\metric{\ad_u v}{w}+\metric{v}{\ad_uw}=0$ for all $u \in \k$, and a naturally reductive decomposition implies
for $u,v,w\in\p$,  $\metric{[u,v]_\p}{w}= \metric{u}{[v,w]_\p}$, where $\p = \k^{\perp}$.

  We want to show that $\ad^*_vv=0$ for any $v \in \p$.
   Let $z \in \g$. 
 Then \begin{align*}
 \metric{\ad^*_vv}{z}&=\metric{\ad^*_vv}{z_\k}+\metric{\ad^*_vv}{z_\p}=\metric{v}{[v,z_\k]}+\metric{v}{[v,z_\p]} \\
 &=
 \metric{v}{[v,z_\k]}+\metric{[v,v]}{z_\p}= \metric{v}{[v,z_\k]}=- \metric{v}{\ad_{z_\k}v},
 \end{align*} 
 where we have used $$\metric{v}{[v,z_\p]}=\metric{v}{[v,z_\p]_\p}=\metric{[v,v]_\p}{z_\p}=\metric{[v,v]}{z_\p}=0,$$
  by orthogonality and the fact $[v,v]=0$.
 
 Moreover $ \metric{v}{\ad_{z_\k}v}=-\metric{\ad_{z_\k}v}{v}=- \metric{v}{\ad_{z_\k}v}$, so $\metric{v}{\ad_{z_\k}v}=0$. In fact this only requires the inner product on $\g$ to be $\ad_\k\p$-invariant, that is $\metric{\ad_u v}{w}+\metric{v}{\ad_uw}=0$ for all $u \in \k, v,w \in \p$.

Finally let us consider the term $\ad^*_{v_\p}v_\k$. 
Note let $S_i, A_j$ basis of $\p,\k$ respectively, and $T_k$ arbitrary basis element of $\g$. 
Then $\metric{\ad^*_{S_i}A_j}{T_k} = \metric{A_j}{[S_i,T_k]}$. This is zero if $T_k \in \k$, since $[\p,\k] \subset \p$ and $\p=\k^{\perp}$.
Moreover 
$\metric{\ad^*_{S_i}A_j}{T_k}= \metric{A_j}{c^r_{ik}T_r}= c^{r}_{ik}g_{jr}=c_{jik}$
where $j \in \k$, $i \in \p$.
On the other hand
$\metric{[S_i,A_j]}{T_k}$ is again zero if $T_k \in \k$; and $\metric{[S_i,A_j]}{T_k}= c^r_{ij}g_{rk}=c_{kij} = -c_{kji}=c_{ijk}$ for $i \in \p$, $j\in \k$, where in the last step we used $\ad_\k$-invariance.
It doesn't seem to be possible in general to relate $c_{jik}$ with $c_{ijk}$.
\subsection{Geodesics Matrix Groups} \label{geo-matrix}

Let us first verify that $\ad^*_SA=[S,A]$, for $S\in \p$, $A\in \k$.
Let $T\in \g$. Then
\begin{align*}
 \metric{\ad^*_SA}T&= \metric{A}{[S,T]}= \tr\big(A^T[S,T]\big)=-\tr\big(A[S,T]\big)\\
 &= -\tr\big(AST\big)+\tr\big(ATS\big)
 = -\tr\big(AST\big)+\tr\big(SAT\big) \\
 &= \tr\big([S,A]T \big)=\metric{[S,A]}{T}
\end{align*}
and since $T$ is arbitrary the result holds. 
Here we have used the property of the trace that $\tr(ABC)=\tr(BCA)$ and the fact that $[S,A]$ is symmetric, 
$[S,A]^T=\big(SA-AS)^T=-AS+SA=[S,A]$.
It follows that the Euler-Arnold equation becomes
 $\dot v = \dot v_\k+ \dot v_\p= -2 [v_\k,v_\p]$. 
 Recall $\ad_u = \tang 1 \Ad (u)= \frac{d}{dt} \Ad(e^{tu})|_{t=0}$ for any $u\in \g$.
 
 For any $u,w \in \g$ the curve $t\mapsto \Ad(e^{tu})(w)$ is tangent to $[u,\Ad(e^{tu})(w)]$ at $t$, since (for matrix groups)
 $\Ad(e^{tu})w=e^{tu}we^{-tu}$ and differentiating gives $ue^{tu}we^{-tu}-e^{tu}we^{-tu}u=u\cdot\Ad(e^{tu})w-\Ad(e^{tu})w\cdot u=[u,\Ad(e^{tu})w]$;
 in particular at $t=0$ this is $[u,w]$. 
 Moreover note that $\Ad \K (\p) \subset \p$ (i.e., the decomposition of $\g$ is reductive).
 Hence if $v(0)=v_\k(0)+v_\p(0)$ and using $\Ad(1)(u)=u$, we find that
 the curve $$v(t) =v_\k(0)+ \Ad e^{-2 v_\k(0)t}(v_\p(0))$$ satisfies the Euler-Arnold equation \big(i.e. $v_\k(t)= v_\k(0)$ and $v_\p(t) =\Ad e^{-2 v_\k(0)t}(v_\p(0))$\big).
 
The reconstruction equation then implies~\cite{bryant1995introduction}
$$q(t)= q(0)e^{(v_\p-v_\k)t}e^{2v_\k t}.$$
 Note we can also write 
 $$v_\p - v_\k = \half (v+v^T)-\half  (v-v^T)=v^T, \qquad 2v_\k = v-v^T.$$

\subsection{Symmetries of Hamiltonian} \label{sym-Ham}

For $v\in \g$, the Nother's current $J$ associated to the lifted right action of $\K$ on $T\G$ is $J \circ \legtr( v) = \metric{v}{\cdot}_1 \in \k^*$ (see cororally 4.2.11 \cite{abraham1978foundations}). 
Since $H$ is invariant under the action, by theorem 4.2.2 we know that $J$ is conserved along the Hamiltonian flow. 
In particular if $v(0)\in \p= \k^{\perp}$, i.e., $J\circ \legtr(v_0)=0$, then 
$J\circ \legtr(v(t)) \in \p$ for all time.

\subsection{Symplectic Integrators and Shadow Hamiltonian} \label{Symplectic-integrator}

That the flow of a Hamiltonian vector field preserve the symplectic structure is proposition 3.3.4 \cite{abraham1978foundations}.

In order to simulate Hamilton's equations we need to build flows that approximate the Hamiltonian flow $\dot \gamma = \{ \gamma,H\}=\hvf H(\gamma)$.
Given an initial point $\gamma(0) \in \G \times \g$, the solution to Hamilton's equations satisfy
$\gamma(t) = \Phi^{\hvf H}_{\gamma(0)}(t)\defn
 \Phi^{\hvf H}(t) \big(\gamma(0)\big)$.
We say an approximate flow $\Psi $ is said to be a symplectic integrator if it preserves the symplectic structure $\Psi^* \ftfl = \ftfl$ 
(technically we only need integrators that preserve the volume $\ftfl^n$).
Since Hamiltonian vector fields preserve the symplectic structure, a simple way to build integrators of $H$ consists of alternating between the flows of simpler Hamiltonians.
For example the leapfrog integrator is defined as 
$\Psi_{\gamma(0)}(h) \defn  \Phi^{\hvf V}( h/2) \circ \Phi^{\hvf T}(h) \circ \Phi^{\hvf V}_{\gamma(0)}(h/2)$.
This integrator has the important property to also be reversible.
\subsection{Lie Algebra} \label{lie-algebra}

The tangent space $T_1\G$ may be identified with the vector space of left-invariant vector fields $\g_L$, where we associate to each $\xi \in T_1\G$ the left-invariant vector field $X: g \mapsto \tang{1}L_g \xi$, where left-invariance means $(L_g)_*X=X$ and $(L_g)_*$ is the push-forward. The lie algebra $\g$ is $T_1\G$ equipped with the Lie bracket induced by the commutator of vector fields and the isomorphism $T_1\G \cong \g_L$.

\subsection{Horizontal Space} \label{horizontal}

Recall that a principal bundle connection on $\G \to \G/\K$ is simply a smooth choice of horizontal spaces $H_g\subset T_g\G$, s.t. $T\G = ver_g \oplus H_g$ and $ \tang g R_k \big(H_{g}\G\big)=H_{gk}\G$.
In particular we set $\p \defn H_1$.

If $v_1, v_2 \in \p$ and $\tang{1}\pi(v_1)=\tang{1}\pi(v_2)$, then $\tang{1}\pi(v_1-v_2)=0$ 
so $v_1-v_2 \in \k$ and thus $v_1-v_2=0$, since $\k$ and $\p$ are complementary subspaces. Hence $v_1=v_2$.

\subsection{Force Integrators} \label{Force-proof}

Using \refeq{hvf-f}, and the fact
$\{r,h\} = \ftfl(\hvf r, \hvf h) = \hvf h[r]$, we find 
\begin{align*}
\{r,h\} &=  g^{jk} \xi_k(h) e_j(r) + \Big( g^{ik}g^{lj} v^r c_{rji} \xi_l(h)-g^{jk}e_j(h) \Big) \xi_k (r) \\
&= g^{jk}\big( \xi_k(h)e_j(r)-e_j(h)\xi_k(r)\big)+v^rc_r^{lk}\xi_l(h)\xi_k(r)
\end{align*}
where $c_r^{lk} \defn g^{ik}g^{lj} c_{rji}$.
In particular 
$$\{V,T\}= g^{jk}\xi_k(T)e_j(V)=g^{jk}g_{kr}v^re_j(V)=v^re_r(V)$$
and $\{V,\{V,T\}\}= g^{jk} \xi_k \big( v^r e_r(V)\big) e_j(V)= g^{jk} e_k(V) e_j(V)$.
From \refeq{hvf-f} we find 
$$\hvf{\{V,\{V,T\}\}}= -g^{jk}g^{ls} e_j\big( e_l(V) e_s(V)\big)\xi_k = -2g^{jk}g^{ls} e_l(V) e_je_s(V)  \xi_k $$
where we use Leibniz Rule.

To see that $ -2g^{jk}g^{ls} e_l(V) e_je_s(V)  \xi_k\in \p$, we first observe that 
$e_je_s(V)= [e_j,e_s](V)-e_se_j(V)$.
Moreover we can restrict to $s\in \p$ since otherwise $e_s(V)=0$ (note here the basis fields $(e_i)$ are generated from an orthogonal decomposition $\g = \k \oplus \p$).
Then $e_j(V)=0$ for all $j\in \k$, and using the properties of the stucture constant we have
$[e_j,e_s]= c^r_{js}e_r$. 
We have $r \in \p$ (again we have $e_r(V)=0$), and note 
$$g^{jk}g^{ls} e_l(V) [e_j,e_s](V)  \xi_k
=g^{jk}g^{ls} e_l(V) c^r_{js}e_r(V)  \xi_k
= g^{jk}g^{ls} g^{ra} e_l(V) c_{ajs}e_r(V)  \xi_k$$
Hence, since $\p = \k^{\perp}$, $k\in \k$ implies $j\in \k$, implies $c_{ajs}=-c_{sja}$ and since 
$g^{ls} g^{ra} e_l(V)e_r(V)$ is symmetric in $a,s$ the term sums to $0$. 
We have thus shown that $\hvf{\{V,\{V,T\}\}} \in \p$.

\end{document}